\documentclass[12pt,a4paper,fleqn]{article}

\usepackage{amsmath,amsfonts,amsthm,amssymb,color}
\usepackage{latexsym}
\usepackage[latin1]{inputenc}
\usepackage{amsmath}
\usepackage{amsfonts}
\usepackage{oldstyle}
\usepackage{stmaryrd}
\usepackage{color}
\usepackage{delarray}

\usepackage{eufrak}

\newtheorem{prop}{Proposition}[section]

\newtheorem{lemma}[prop]{Lemma}

\newtheorem{theorem}[prop]{Theorem}

\renewcommand{\geq}{\geqslant}
\def\leq{\leqslant}

\newcommand{\R}{\mathbb{R}}

\newcommand{\W}{\mathbb{W}}

\newcommand{\iou}{\int_{0}^{1}}

\newcommand{\intR}{\int_{\R}}

\newcommand{\lp}{\left(}
\newcommand{\rp}{\right)}
\newcommand{\lc}{\left[}
\newcommand{\rc}{\right]}

\newcommand{\lln}{\left|}
\newcommand{\rrn}{\right|}
\newcommand{\lla}{\left\langle}
\newcommand{\rra}{\right\rangle}

\def\HH{\EuFrak H}
\def\W{{{\HH}^{\otimes 2}}}

\def\1{{\mathbf{1}}}

\def\1{{\mathbf{1}}}
\def\0.5{{\frac{1}{2}}}

\begin{document}

\begin{center}
{\large{\bf 
Total variation distance between two\\ 
double Wiener-It\^o integrals
}}\\~\\
Rola Zintout\footnote{Universit\'e de Lorraine, Institut \'Elie Cartan de Lorraine, UMR 7502, Vandoeuvre-l\`es-Nancy, F-54506, France and CNRS, Institut \'Elie Cartan de Lorraine, UMR 7502, Vandoeuvre-l\`es-Nancy, F-54506, France; {\tt rola.zintout@hotmail.com}}\\
{\it Universit\'e de Lorraine}\\~\\
\end{center}

{\small
\noindent
{\bf Abstract:} Using an approach recently developed by Nourdin and Poly \cite{NP}, we improve the rate in an inequality for the total variation distance between two
double Wiener-It\^o integrals originally due to Davydov and Martynova \cite{DM}.
An application to the rate of convergence of a functional of a correlated two-dimensional fractional Brownian motion
towards the Rosenblatt random variable is then given, following a previous study by Maejima and Tudor \cite{MT}.\\



\section{Introduction}

Suppose that $X=\{X(h),\,h\in \EuFrak H\}$ is an isonormal Gaussian process on
a real separable infinite-dimensional Hilbert space $\HH$. For any integer $p\geq 1$, let $%
\EuFrak H^{\otimes p}$ be the $p$th tensor product of $\EuFrak H$. Also, denote
by $\EuFrak H^{\odot p}$ the $p$th symmetric tensor product.

The following statement is due to Davydov and Martynova \cite{DM},
see also \cite[Theorem 4.4]{NP}.
\begin{theorem}\label{dm-thm}
Fix an integer $p\geq 2$, and let $(f_n)$ be a sequence of $\HH^{\odot p}$
 that converges to $f_\infty$ in $\HH^{\otimes p}$.
 Assume moreover that $f_\infty$ is not identically zero.
 let $I_p(f_n)$, $n\in\mathbb{N}\cup\{\infty\}$, denote the $p$th Wiener-It\^o integral of $f_n$ with respect to $X$.
Then, there exists $c>0$ such that, for all $n$,
\begin{equation}\label{dm}
\sup_{C\in\mathcal{B}(\R)}\big|
P(I_p(f_n)\in C)-P(I_p(f_\infty)\in C)
\big|\leq c\,\| f_n - f_\infty\|^{1/p}_{\HH^{\otimes p}},
\end{equation}
where $\mathcal{B}(\R)$ stands for the set of Borelian sets of $\R$.
\end{theorem}
In this paper, $p=2$ and the inequality (\ref{dm}) becomes: 
\begin{equation}\label{DM}
\sup_{C\in\mathcal{B}(\R)}\big|
P(I_2(f_n)\in C)-P(I_2(f_\infty)\in C)
\big|\leq c\,\sqrt{\| f_n - f_\infty\|_{\HH^{\otimes 2}}}.
\end{equation}

To each $f_\infty\in\HH^{\odot 2}$, one may associate the following Hilbert-Schmidt operator:
\begin{equation}\label{Af}
A_{f_\infty}:\HH\to \HH,\quad g\mapsto \langle f_\infty,g\rangle_\HH.
\end{equation}
Let $\lambda_{\infty,k}$, $k\geq 1$, indicate the eigenvalues of $A_{f_\infty}$.
In many situations of interest (see below for an explicit example), it happens that the following  property, that  we label for further use,
is satisfied for $f_\infty$:
\begin{equation}\label{H}
\mbox{the cardinality of $\{k:\,\lambda_{\infty,k}\neq 0\}$ is at least 5. }
\end{equation}
The aim of this paper is  to take advantage of (\ref{H}) in order to improve (\ref{DM}) by a factor 2.
More precisely, relying on an approach recently developed by Nourdin and Poly in \cite{NP}, we shall prove the following result, compare with (\ref{DM}):
\begin{theorem}\label{main}
Let $f_\infty$ be an element of $\HH^{\odot 2}$ satisfying (\ref{H}) (in particular, $f_\infty$ is not identically zero).
Let $(f_n)$ be a sequence of $\HH^{\odot 2}$ that converges to $f_\infty$ in $\HH^{\otimes 2}$.
Then, there exists $c>0$ (depending only on $f_\infty$) such that, for all $n$,
\begin{equation}\label{ineq}
\sup_{C\in\mathcal{B}(\R)}\big|
P(I_2(f_n)\in C)-P(I_2(f_\infty)\in C)
\big|\leq c\| f_n - f_\infty\|_{\HH^{\otimes 2}}.
\end{equation}
\end{theorem}

In some sense, the inequality (\ref{ineq}) appears to be optimal.
Indeed, consider $F_\infty=I_2(f_\infty)$ with $f_\infty$ satisfying
(\ref{H}) and set $F_n=I_2(f_n)$ with $f_n=(1+c_n)f_\infty$, where
$(c_n)$ is a sequence of nonzero real numbers converging to zero.
Let $\phi_\infty$ (resp. $\phi_n$) denote
the density of $F_\infty$ (resp. $F_n$), which exists thanks to Shigekawa's theorem (see \cite{S}). Assume furthermore that $\phi_\infty$ is differentiable and is such that
$0<\int_\R |x\phi_\infty'(x)+\phi_\infty(x)|dx<\infty$.
According to Scheff\'e's theorem, one has
\[
\sup_{C\in\mathcal{B}(\R)}\big|
P(I_2(f_n)\in C)-P(I_2(f_\infty)\in C)\big|
=\frac12\int_\R |\phi_n(x))-\phi_\infty(x)|dx.
\]
We deduce, after some easy calculations, that
\[
\sup_{C\in\mathcal{B}(\R)}\big|
P(I_2(f_n)\in C)-P(I_2(f_\infty)\in C)\big|
\sim_{n\to \infty}\frac12|c_n|\int_\R |x\phi_\infty'(x)+\phi_\infty(x)|dx.
\]
On the other hand, $\| f_n - f_\infty\|_{\HH^{\otimes 2}}=|c_n|\,\,\| f_\infty\|_{\HH^{\otimes 2}}$. Thus,
\[
\sup_{C\in\mathcal{B}(\R)}\big|
P(I_2(f_n)\in C)-P(I_2(f_\infty)\in C)\big|
\sim_{n\to \infty} c\,\| f_n - f_\infty\|_{\HH^{\otimes 2}},
\]
with $c=\int_\R |x\phi_\infty'(x)+\phi_\infty(x)|dx/
(2\| f_\infty\|_{\HH^{\otimes 2}})$.

To illustrate the use of Theorem \ref{main} in a concrete situation, we consider the following example taken from Maejima and Tudor \cite{MT}.
Let $B^{H_1}$, $B^{H_2}$ be two fractional Brownian motions with Hurst parameters $H_1,H_2 \in (0,1)$, respectively.
We assume that both $H_1$ and $H_2$ are strictly bigger than $\frac12$. We further assume that
the two fractional Brownian motions $B^{H_1}$ and $B^{H_2}$ can be expressed as Wiener integrals with respect to the {\it same} two-sided Brownian motion $W$, meaning in particular that $B^{H_1}$ and $B^{H_2}$ are {\it not} independent. Precisely, we set
\begin{eqnarray}\label{fbm1}
B^{H_1}_t &=&c(H_1)\int_\R dW_y\int_0^t (u-y)_+^{H_1-\frac32}du,\quad t\geq 0, \\
B^{H_2}_t &=&c(H_2)\int_\R dW_y\int_0^t (u-y)_+^{H_2-\frac32}du, \quad t\geq 0, 
\label{fbm2}
\end{eqnarray}
where the constants $c(H_1)$ and $c(H_2)$ are chosen so that $E[(B^{H_1}_1)^2]=E[(B^{H_2}_1)^2]=1$.
Define
\begin{equation}\label{Zn}
Z_n=n^{1-H_1-H_2}\sum_{k=0}^{n-1}\lc \frac{(B_{\frac{k+1}{n}}^{H_1}-B_{\frac{k}{n}}^{H_1})
(B_{\frac{k+1}{n}}^{H_2}-B_{\frac{k}{n}}^{H_2})}{E\lc (B_{\frac{k+1}{n}}^{H_1}-B_{\frac{k}{n}}^{H_1})
(B_{\frac{k+1}{n}}^{H_2}-B_{\frac{k}{n}}^{H_2})\rc}-1\rc.
\end{equation}
When $H_1=H_2=H$, observe that (\ref{Zn}) is related to the quadratic
variation of $B^H$.
In \cite{MT}, the following extension of a classical result by Taqqu \cite{Taqqu} is shown:

\begin{prop}\label{maejima-tudor}
Assume that $H_1>\frac12$, $H_2>\frac12$ and $H_1+H_2>\frac32$.
Then, $Z_n$ converges as $n\to\infty$ in $L^2(\Omega)$ to
 the non-symmetric Rosenblatt random variable
 $Z_\infty$, given by
\begin{equation}\label{Zinfty}
Z_\infty = b(H_1,H_2)\int_{\R^2} dW_{x}dW_{y}\,\int_0^1 (s-x)_+^{H_1-3/2}(s-y)_+^{H_2-3/2}ds.
\end{equation}
Here $b(H_1,H_2)$ is a normalizing explicit constant
whose precise value does not matter in the sequel.
\end{prop}

In the present paper, by relying on (\ref{ineq}) we are able to associate an explicit rate to the convergence
$Z_n\overset{L^2}{\to}Z_\infty$ of Proposition \ref{maejima-tudor}, namely,
\begin{equation}\label{rate}
\sup_{C\in\mathcal{B}(\R)}\big|
P(Z_n\in C)-P(Z_\infty\in C)
\big|=O(n^{\frac32-H_1-H_2}).
\end{equation}
When $H_1=H_2=H$, the rate $\frac32-2H$ we have obtained in (\ref{rate}) is better (by a power 2) than the one computed by Breton and Nourdin in \cite{BN}, precisely because our inequality (\ref{ineq}) improves the inequality (\ref{DM})
of Davydov and Martynova by a power 2.\\

The rest of the paper is organized as follows. Section 2 contains some preliminary material on Malliavin calculus. In Section 3 we prove Theorem \ref{main}. Finally, Section 4 contains our proof of (\ref{rate}).
\setcounter{equation}{0}
\section{Preliminaries}

Let $\EuFrak H$ be a real separable infinite-dimensional Hilbert space. For any integer $p\geq 1$, let $%
\EuFrak H^{\otimes p}$ be the $p$th tensor product of $\EuFrak H$. Also, we denote
by $\EuFrak H^{\odot p}$ the $p$th symmetric tensor product.

Suppose that $X=\{X(h),\,h\in \EuFrak H\}$ is an isonormal Gaussian process on
$\EuFrak H$, defined on some probability space $(\Omega ,\mathcal{F},P)$.
Assume from now on that $\mathcal{F}$ is
generated by $X$.
For every integer $p\geq 1$, let $\mathcal{H}_{p}$ be the $p$th Wiener chaos of $X$,
that is, the closed linear subspace of $L^{2}(\Omega)$
generated by the random variables $\{H_{p}(X(h)),h\in \EuFrak H,\left\|
h\right\| _{\EuFrak H}=1\}$, where $H_{p}$ is the $p$th Hermite polynomial defined by
\[
H_p(x)=\frac{(-1)^p}{p!}e^{x^2/2}\frac{d^p}{dx^p}\big(e^{-x^2/2}\big).
\]
We denote by $\mathcal{H}_{0}$ the space of constant random variables. For
any $p\geq 1$, the mapping $I_{p}(h^{\otimes p})=p!H_{p}(X(h))$,
$h\in \EuFrak H$, $\left\|
h\right\| _{\EuFrak H}=1$, provides a
linear isometry between $\EuFrak H^{\odot p}$
(equipped with the modified norm $\sqrt{p!}\left\| \cdot \right\| _{\EuFrak %
H^{\otimes p}}$) and $\mathcal{H}_{p}$ (equipped with the $L^2(\Omega)$ norm). 
We call $I_p(f)$ the $p$th multiple Wiener-It\^o integral of kernel $f$.
For $p=0$, by convention $%
\mathcal{H}_{0}=\mathbb{R}$, and $I_{0}$ is the identity map.
In particular, when $f,g\in \EuFrak H^{\odot p}$, observe that
\begin{equation}\label{isometry}
E\lc \lp I_p(f)-I_p(g)\rp^2\rc = p!\left\|f-g\right\|^2 _{\EuFrak H^{\otimes p}}.
\end{equation}
It is well-known (Wiener chaos expansion) that $L^{2}(\Omega)$
can be decomposed into the infinite orthogonal sum of the spaces $\mathcal{H}%
_{p}$. That is, any square integrable random variable $F\in L^{2}(\Omega)$ admits the following chaotic expansion:
\begin{equation}
F=\sum_{p=0}^{\infty }I_{p}(f_{p}),  \label{E}
\end{equation}%
where $f_{0}=E[F]$, and the $f_{p}\in \EuFrak H^{\odot p}$, $p\geq 1$, are
uniquely determined by $F$. For every $p\geq 0$, we denote by $J_{p}$ the
orthogonal projection operator on the $p$th Wiener chaos. In particular, if $%
F\in L^{2}(\Omega)$ is as in (\ref{E}), then $%
J_{p}F=I_{p}(f_{p})$ for every $p\geq 0$.

Let us now introduce some basic elements of the Malliavin calculus with respect
to the isonormal Gaussian process $X$. We refer the reader to Nourdin and Peccati \cite{NouPecBook} or Nualart \cite%
{nualartbook} for a more detailed presentation of these notions. Let $\mathcal{S}$
be the set of all smooth and cylindrical random variables of
the form
\begin{equation}
F=g\left( X(\phi _{1}),\ldots ,X(\phi _{n})\right) ,  \label{v3}
\end{equation}%
where $n\geq 1$, $g:\mathbb{R}^{n}\rightarrow \mathbb{R}$ is an infinitely
differentiable function with compact support, and $\phi _{i}\in \EuFrak H$.
The Malliavin derivative of $F$ with respect to $X$ is the element of $%
L^{2}(\Omega ,\EuFrak H)$ defined as
\begin{equation*}
DF\;=\;\sum_{i=1}^{n}\frac{\partial g}{\partial x_{i}}\left( X(\phi
_{1}),\ldots ,X(\phi _{n})\right) \phi _{i}.
\end{equation*}
By iteration, one can
define the $k$th derivative $D^{k}F$ for every $k\geq 2$, which is an element of $L^{2}(\Omega ,%
\EuFrak H^{\odot k})$.

For $k\geq 1$ and $p\geq 1$, ${\mathbb{D}}^{k,p}$ denotes the closure of $%
\mathcal{S}$ with respect to the norm $\Vert \cdot \Vert_{\mathbb{D}^{k,p}}$, defined by
the relation
\begin{equation*}
\Vert F\Vert _{\mathbb{D}^{k,p}}^{p}\;=\;E\left[ |F|^{p}\right] +\sum_{i=1}^{k}E\left(
\Vert D^{i}F\Vert _{\EuFrak H^{\otimes i}}^{p}\right) .
\end{equation*}
The Malliavin derivative $D$ verifies the following chain rule. If $%
\varphi :\mathbb{R}^{n}\rightarrow \mathbb{R}$ is continuously
differentiable with bounded partial derivatives and if $F=(F_{1},\ldots
,F_{n})$ is a vector of elements of ${\mathbb{D}}^{1,2}$, then $\varphi
(F)\in {\mathbb{D}}^{1,2}$ and
\begin{equation}\label{chain}
D\varphi (F)=\sum_{i=1}^{n}\frac{\partial \varphi }{\partial x_{i}}%
(F)DF_{i}.
\end{equation}
Observe that (\ref{chain}) still holds when $\varphi$ is Lipschitz and the law of $F$ has a density with respect to the Lebesgue measure on $\R^n$ (see, e.g., Proposition 1.2.3 in \cite{nualartbook}).

We denote by $\delta $ the adjoint of the operator $D$, also called the 
divergence operator.
A random element $u\in L^{2}(\Omega ,\EuFrak %
H)$ belongs to the domain of $\delta $, noted $\mathrm{Dom}\delta $, if and
only if it verifies
\begin{equation*}
\big|E\big(\langle DF,u\rangle _{\EuFrak H}\big)\big|\leq c_{u}\,\sqrt{E(F^2)}
\end{equation*}%
for any $F\in \mathbb{D}^{1,2}$, where $c_{u}$ is a constant depending only
on $u$. If $u\in \mathrm{Dom}\delta $, then the random variable $\delta (u)$
is defined by the duality relationship:
\begin{equation}
E(F\delta (u))=E\big(\langle DF,u\rangle _{\EuFrak H}\big),  \label{ipp}
\end{equation}
which holds for every $F\in {\mathbb{D}}^{1,2}$. We will also make use of the following relationships, valid
for 
$F\in \mathbb{D}^{1,2}$ and $u\in \mathrm{Dom}\delta $ such that $Fu\in L^{2}(\Omega ,\EuFrak %
H)$:
\begin{eqnarray}
F\delta (u)&=&\delta \lp Fu\rp +\langle DF,u\rangle _{\EuFrak H} \label{ipp2}\\
E\lp \delta(u)^2\rp&=&E\left\| Du \right\| ^2_{\EuFrak H^{\otimes 2}} + E \left\| u \right\| ^2_{\EuFrak H}. \label{varianceskorohod}
\end{eqnarray}
The operator $L$ is defined on the Wiener chaos expansion as
\begin{equation*}
L=\sum_{q=0}^{\infty }-qJ_{q},
\end{equation*}
and is called the infinitesimal generator of the Ornstein-Uhlenbeck
semigroup. The domain of this operator in $L^{2}(\Omega)$ is the set%
\begin{equation*}
\mathrm{Dom}L=\{F\in L^{2}(\Omega ):\sum_{q=1}^{\infty }q^{2}\left\|
J_{q}F\right\| _{L^{2}(\Omega )}^{2}<\infty \}=\mathbb{D}^{2,2}\text{.}
\end{equation*}%
There is an important relationship between the operators $D$, $\delta $ and $L$. A random variable $F$ belongs to the
domain of $L$ if and only if $F\in \mathrm{Dom}\left( \delta D\right) $
(i.e. $F\in {\mathbb{D}}^{1,2}$ and $DF\in \mathrm{Dom}\delta $), and in
this case
\begin{equation}
\delta DF=-LF.  \label{LF}
\end{equation}

If $\EuFrak H=%
L^{2}(A,\mathcal{A},\mu )$ (with $\mu $ non-atomic), then the
derivative of a random variable $F$ as in (\ref{E}) can be identified with
the element of $L^{2}(A\times \Omega )$ given by
\begin{equation}
D_{a}F=\sum_{q=1}^{\infty }qI_{q-1}\left( f_{q}(\cdot ,a)\right) ,\quad a\in
A.  \label{dtf}
\end{equation}
At this stage, we observe that an easy calculation leads to the following
identity for $F=I_p(f)$ and $G=I_p(g)$ (with $f,g\in \EuFrak H^{\odot p}$), that we label for further use:
\begin{equation}\label{square-norm-malliavin}
E\lp\left\|DF-DG\right\|_{\HH}^2\rp=pp!\left\|f-g\right\|^2_{\EuFrak H^{\otimes p}}.
\end{equation}

Finally, the following lemma will play a crucial role in our forthcoming calculations.
\begin{lemma}\label{lem}
Let $F_\infty =I_2(f_\infty)$, with $f_\infty\in\HH^{\odot 2}$ satisfying (\ref{H}).
Then, for all $r\geq 1$, we have
\begin{eqnarray}\label{hyper}
E[|F_\infty|^{2r}]<\infty ,\quad E\lc \| DF_\infty\|_{\HH}^{2r}\rc<\infty,
\end{eqnarray} 
as well as
\begin{eqnarray}\label{inverse}
E\lc \frac{1}{\left\| DF_\infty\right\|_{\HH}^{9/2}}\rc<\infty.
\end{eqnarray} 
\end{lemma}
\noindent
{\it Proof}.
The proof of (\ref{hyper}) is classical and follows directly from the hypercontractivity property of multiple Wiener-It\^o integrals. So, let us only focus on (\ref{inverse}).
Let $e_{k}$, $k\geq 1$, be the eigenvectors associated to the eigenvalues $\lambda_{f_\infty,k}$ of $A_{f_\infty}$, see (\ref{Af}).
Observe that they form an orthonormal system in $\HH$ and that $f_\infty$ may be expanded as
\begin{equation}\label{finfini}
f_\infty=\sum_{k=1}^\infty \lambda_{f_\infty,k}\,e_k\otimes e_k,
\end{equation}
implying in turn that
\[
F_\infty=I_2(f_\infty) = \sum_{k=1}^\infty \lambda_{f_\infty,k} \,(X(e_k)^2-1).
\]
We have
\begin{eqnarray*}
E\lc \frac{1}{\left\| DF_\infty\right\|_{\HH}^{9/2}}\rc&=&
\int_0^\infty P\left(
\frac{1}{\left\| DF_\infty\right\|_{\HH}^{9/2}}\geq x
\right)
dx\\
&=&\int_0^1 P\left(
\frac{1}{\left\| DF_\infty\right\|_{\HH}^{9/2}}\geq x
\right)
dx+\int_1^\infty P\left(
\frac{1}{\left\| DF_\infty\right\|_{\HH}^{9/2}}\geq x
\right)
dx\\
&\leq&
1+\int_1^\infty P\left(
\left\| DF_\infty\right\|_{\HH}^{2}\leq x^{-4/9}
\right)
dx\\
&\leq&
1+\frac94\int_0^1 P\left(
\| DF_\infty\|_{\HH}^{2}\leq u
\right)\frac{du}{u^{13/4}}.
\end{eqnarray*}
To achieve the desired conclusion (\ref{inverse}), let us check that
\begin{equation}\label{toshow}
P\left(
\| DF_\infty\|_{\HH}^{2}\leq u
\right)=O(u^{5/2})\quad\mbox{as $u\downarrow 0$}.
\end{equation}
An immediate calculation leads to
\begin{equation}\label{DF}
\|DF_\infty\|^2_{\HH} = 4 \sum_{k=1}^\infty \lambda_{f_\infty,k}^2\,X(e_k)^2,
\end{equation}
where the $X(e_k)$ are independent $N(0,1)$ random variables.
Therefore, for any $u>0$,
\begin{eqnarray*}
P\left(
\| DF_\infty\|_{\HH}^{2}\leq u
\right)&\leq&
 P\left(\bigcap_{i=1}^5\{4\lambda_{f_\infty,i}^2X(e_i)^2\leq u\}\right)
=\prod_{i=1}^5 P\left(|X(e_i)|\leq \frac{\sqrt{u}}{2|\lambda_{f_\infty,i}|}\right)\\
&\leq&\frac{u^{5/2}}{(2\pi)^{5/2}\prod_{i=1}^5|\lambda_{f_\infty,i}|}
\end{eqnarray*}
and (\ref{toshow}) is checked, thus concluding the proof.
\qed

\section{Proof of Theorem \ref{main}}

Let $f_\infty\in\HH^{\odot 2}$ satisfying (\ref{H}) (in particular,
$f_\infty$ is not identically zero).
Let $(f_n)$ be a sequence of $\HH^{\odot 2}$ that converges to $f_\infty$ in $\HH^{\otimes 2}$.
Write $F_n=I_2(f_n)$ and $F_\infty=I_2(f_\infty)$.
Our aim in this Section 3 is to show that there exists $c>0$ (depending only on $f_\infty$) such that, for all Borelian set $C$ and all $n$,
\begin{equation}\label{main2}
\big|
P(F_n\in C)-P(F_\infty\in C)
\big|\leq c\| f_n - f_\infty\|_{\HH^{\otimes 2}}.
\end{equation}

First of all, relying on the Lebesgue's monotone convergence theorem, we notice that it is not a loss of generality to assume that the Borel set $C$ is bounded in (\ref{main2}).

Now, we split the proof of Theorem \ref{main} into several steps and we stress that, in what follows, the constant $c$ shall denote a generic constant only depending on $f_\infty$ (not on $n$ !) and whose value may change
from one line to another. 

\bigskip

\underline{\it Step 1}. Thanks to (\ref{DF}), we have
$
\|DF_\infty\|^2_{\HH}\geq 4 \lambda_{f_\infty,k}^2\,X(e_k)^2
$
for some $k$ with $\lambda_{f_\infty,k}\neq 0$ (assumption (\ref{H})). Since $X(e_k)\neq 0$ a.s., one has that 
$\|DF_\infty\|_{\HH}>0$ a.s.
As a result, one can write
\begin{eqnarray} \label{distance}
\lln P\lp F_n\in C\rp -P\lp F_\infty \in C\rp \rrn = 
\lln E\lc \lp {1}_{F_n\in C}-{1}_{F_\infty \in C}\rp \frac{\| DF_{\infty} \|^2_{\HH}}
{\| DF_{\infty} \|^2_{\HH}}
\rc \rrn .
\end{eqnarray}
The chain rule for Lipschitz function (for $n$ large enough, note that $F_n$ has a density with respect to the Lebesgue measure by Shigekawa theorem \cite{S}) leads to
\[
D(\int_{-\infty}^{F_n}1_C(x)dx)=1_C(F_n)DF_n\quad\mbox{and}\quad
D(\int_{-\infty}^{F_\infty}1_C(x)dx)=1_C(F_\infty)DF_\infty.
\]
We then have
\begin{equation}
\label{eq3}
\big|
P(F_n\in C)-P(F_\infty\in C)
\big|\leq |A_n|+|B_n|,
\end{equation}
with
\begin{eqnarray}
A_n&=&E
\left[
\frac{\left\langle D\lp \int_{F_\infty}^{F_n}1_C(x)dx \rp,DF_\infty \right\rangle_{\HH}}
{\|DF_\infty\|^2_{\HH}}\rc \label{A}\\
B_n&=&E\lc\frac{1_C(F_n)\lla D\lp F_\infty - F_n\rp ,DF_\infty\rra_{\HH}}{\|DF_\infty\|^2_{\HH}}\rc\label{B}.
\end{eqnarray}

\bigskip

\underline{\it Step 2} (a bound for $B_n$).
Using Cauchy-Schwarz inequality twice, one obtains
\begin{eqnarray*}\label{B2}
|B_n|&\leq&E\lc\frac{\| D(F_\infty - F_n)\|_\HH}{\|DF_\infty\|_\HH}\rc 
\leq \sqrt{E\| D(F_\infty - F_n)\|_\HH^2}\sqrt{E\lc\frac{1}{\|DF_\infty\|_\HH^2} \rc}.
\end{eqnarray*}
By (\ref{square-norm-malliavin}), one has $E\| D(F_\infty - F_n)\|_\HH^2\leq 4\|f_\infty - f_n\|_{\HH^{\otimes 2}}^2$,  whereas
$E\lc\frac{1}{\|DF_\infty\|_\HH^2} \rc$ is finite by Lemma \ref{lem}.
Thus, 
\begin{equation}\label{zz}
|B_n|\leq c \|f_\infty - f_n\|_\W
\end{equation}
with $c$ only depending on $f_\infty$.

\bigskip

\underline{\it Step 3} (a bound for $A_n$).
Using (\ref{ipp2}), (\ref{LF}) and then Cauchy-Schwarz, one can write
\begin{eqnarray*}
A_n&=&
E\left[
\int_{F_\infty}^{F_n} 1_C(x)dx\,\,\delta\left(
\frac{DF_\infty}{\|DF_\infty\|_\HH^2}
\right)
\right]\\
&=&E\left[
\int_{F_\infty}^{F_n} 1_C(x)dx\,\,
\left\{
\frac{2F_\infty}{\|DF_\infty\|^2_\HH}
-\left\langle DF_\infty,D\frac{1}{\|DF_\infty\|_\HH}\right\rangle_\HH
\right\}
\right]\\
&\leq&
\sqrt{E[(F_n-F_\infty)^2]}\times
\sqrt{
8E\left[
\left(
\frac{F_\infty}{\|DF_\infty\|^2_\HH}
\right)^2
\right]
+2
 E\lp \lla DF_\infty , D \lp \frac{1}{\|DF_\infty\|_\HH^2}\rp  \rra_\HH^2 \rp
}.
\end{eqnarray*}
By Lemma \ref{lem}, it is clear that $E\left[
\left(
\frac{F_\infty}{\|DF_\infty\|^2_\HH}
\right)^2
\right]<\infty$.
On the other hand,
one has
\begin{eqnarray*}
&&E\lp \lla DF_\infty , D \lp \frac{1}{\|DF_\infty\|_\HH^2}\rp  \rra_\HH^2 \rp
=
64\,E\lp \lla DF_\infty , \frac{\lla  f_\infty,DF_\infty\rra_\HH }{\|DF_\infty \|_\HH^4} \rra^2_\HH\rp \\
&=& 64\,E\lp  \frac{\lla f_\infty,DF_\infty\otimes DF_\infty \rra_\W^2}{\|DF_\infty \|^8_\HH}\rp
\leq  64\,\|f_\infty \|_\W^2\,E\lp \frac{1}{\|DF_\infty \|_\HH^4} \rp,
\end{eqnarray*}
which is also finite by Lemma \ref{lem}.
Thus, see also (\ref{isometry}), one has
\begin{equation}\label{zzz}
|A_n|\leq c \|f_\infty - f_n\|_\W,
\end{equation}
with $c$ only depending on $f_\infty$.

\bigskip

\underline{\it Step 4} (conclusion).
Taking into account (\ref{eq3}), (\ref{zz}) and (\ref{zzz}), we obtain that (\ref{main2}) holds true, thus concluding the proof of
Theorem \ref{main}.
\qed

\section{Proof of (\ref{rate})}

To prove (\ref{rate}), we shall apply our Theorem \ref{main}. The
isonormal Gaussian process $X=\{X(h):\,h\in\HH\}$ we consider here is
a two-sided Brownian motion $W=\{W(h):\,h\in L^2(\R)\}$.
We divide the proof of (\ref{rate}) into several steps.

\bigskip

\underline{\it Step 1}. 
Recall from (\ref{Zn}) and (\ref{Zinfty}) the definitions of $Z_n$ and
$Z_\infty$ respectively.
In Maejima and Tudor \cite{MT}, the authors represent $Z_n$ and $Z_\infty$ as  
\[
Z_n=b(H_1,H_2)\times I_2(f_n)\quad\mbox{and}\quad 
Z_\infty=b(H_1,H_2)\times I_2(f_\infty),
\] 
with $b(H_1,H_2)$ a suitable constant and
\begin{eqnarray*}
f_n(x,y)&=& n\sum_{i=0}^{n-1}\int_{i/n}^{(i+1)/n} \int_{i/n}^{(i+1)/n} (s-x)_+^{H_1-3/2}(s-y)_+^{H_2-3/2}ds \\
f_\infty (x,y)&=& \int_0^1 (s-x)_+^{H_1-3/2}(s-y)_+^{H_2-3/2}ds.
\end{eqnarray*}
We have moreover, see indeed \cite[page 180]{MT},
\begin{equation}\label{page180}
\|f_n-f_\infty\|_{L^2(\R^2)} = O(n^{\frac{3}{2}-H_1-H_2})\quad\mbox{as $n\to\infty$}.
\end{equation}

\bigskip

\underline{\it Step 2}. 
Let us check that $f_\infty$ satisfies (\ref{H}). To do so, 
recall from (\ref{finfini}) that $f_\infty$ may be expanded, with $e_{k}$  the eigenvectors associated to $\lambda_{f_\infty,k}$,  as
\begin{equation}\label{ide}
f_\infty (x,y)= \sum_{k\geq 1} \lambda_{\infty,k} e_k(x)e_k(y).
\end{equation}

Let us first show that $e_k$ is bounded on $[0,1]$ when $\lambda_{\infty,k}\neq 0$. 
Indeed, using Cauchy-Schwarz inequality as well as the identity
\[
\int_\R (t-x)_+^{\alpha} (s-x)_+^{\alpha} dx = c_\alpha |t-s|^{2\alpha+1}
\]
valid for any $\alpha>-\frac12$ (with $c_\alpha>0$ a constant depending only on $\alpha$),
one can write
\begin{eqnarray*}
e_k(y)^2&=& \frac{1}{\lambda_k^2} \lp\intR e_k(x)dx \iou ds (s-x)_+^{H_1-3/2} (s-y)_+^{H_2-3/2}\rp^2\\
&\leq& 
\frac{1}{\lambda_k^2} \intR e_k(x)^2dx \times \int_\R dx \lp\iou ds (s-x)_+^{H_1-3/2} (s-y)_+^{H_2-3/2}\rp^2\\
\\
&=&\frac{1}{\lambda_k^2} \intR dx \int_{[0,1]^2} dt ds (s-x)_+^{H_1-3/2} (s-y)_+^{H_2-3/2}
(t-x)_+^{H_1-3/2} (t-y)_+^{H_2-3/2}\\
&=&\frac{c_{H_1}}{\lambda_k^2}\int_{[0,1]^2}dt ds (s-y)_+^{H_2-3/2} (t-y)_+^{H_2-3/2}
|t-s|^{{2H_1}-2},\\
\end{eqnarray*}
with $c_{H_1}$ a constant depending only on $H_1$. Thus, for any $0\leq y \leq 1$,
\begin{eqnarray*}
e_k(y)^2&\leq & \frac{c_{H_1}}{\lambda_k^2}\int_{[y,1]^2}dt ds (s-y)^{H_2-3/2} (t-y)^{H_2-3/2}
|t-s|^{{2H_1}-2}\\
&= & \frac{c_{H_1}}{\lambda_k^2}\int_{[0,1-y]^2}dt ds \quad s^{H_2-3/2} t^{H_2-3/2}|t-s|^{{2H_1}-2}\\
&\leq& \frac{c_{H_1}}{\lambda_k^2}\int_{[0,1]^2}dt ds \quad s^{H_2-3/2} t^{H_2-3/2}|t-s|^{2{H_1}-2}\\
&=& \frac{2c_{H_1}}{\lambda_k^2}\iou dt \, t^{2H_1+2H_2-4}\int_{0}^1 du \, u^{H_2-3/2} (1-u)^{2{H_1}-2}<\infty.
\end{eqnarray*}

Let us now show that $f_\infty$ is {\it not} bounded on $[0,1]^2$. 
If $x,y \in [0,\frac12],$ then
\begin{eqnarray*}
&&\int_0^1 (s-x)_+^{H_1-3/2}(s-y)_+^{H_2-3/2}ds = \int_{x\vee y}^1 (s-x)^{H_1-3/2}(s-y)^{H_2-3/2}ds\\
&\geq &\int_{x\vee y}^1 \frac{ds}{\sqrt{(s-x)(s-y)}}
=\int_{x\vee y}^1\lc \lp s-\frac{x+y}{2}\rp^2 - \lp\frac{x-y}{2}\rp^2 \rc^{-\frac{1}{2}} ds\\
&=& \frac{1}{2}\int_{0}^{(1-x)(1-y)}\frac{du}{\sqrt{u\left(u+(\frac{x-y}{2})^2\right)}}
\geq \frac{1}{2} \int_{0}^{\frac14}\frac{du}{\sqrt{u\left(u+(\frac{x-y}{2})^2\right)}}.
\end{eqnarray*}
Using Fatou's lemma, we conclude that
\begin{eqnarray*}
\liminf_{y\rightarrow x} \int_0^1 (s-x)_+^{H_1-3/2}(s-y)_+^{H_2-3/2}ds&\geq &
\liminf_{y\rightarrow x} \frac{1}{2}\int_{0}^{\frac14}\frac{du}{\sqrt{u\left(u+(\frac{x-y}{2})^2\right)}}\\
&\geq &\frac{1}{2}\int_0^{\frac14}\frac{du}{u} = +\infty.
\end{eqnarray*}

The fact that $f_\infty$ is not bounded together with the fact that $e_k$ is bounded when
$ \lambda_{\infty,k}\neq 0$ imply, thanks to (\ref{ide}), that $f_\infty$ satisfies (\ref{H}).

\bigskip

\underline{\it Step 3} (conclusion). Due to the conclusion of Step 2, the proof of (\ref{rate}) now follows from Theorem \ref{main} and (\ref{page180}).\qed

\bigskip
\noindent
{\bf Acknowledgments}. I would like to thank sincerely my thesis advisor Ivan Nourdin, who led all the directions of my work.
Also, I am particularly grateful to two anonymous referees for a careful reading and a number of helpful suggestions, which led to significant improvements in the presentation of my results.

\end{document}